\documentclass[12pt, a4paper]{amsart}

\usepackage[hmargin=30mm, vmargin=25mm, includefoot, twoside]{geometry}
\usepackage[bookmarksopen=true]{hyperref}

\usepackage{amsfonts,amssymb,verbatim}
\usepackage{latexsym}
\usepackage{mathrsfs}
\usepackage{stmaryrd}
\usepackage{xspace}
\usepackage{enumerate, paralist}
\usepackage{graphicx}
\usepackage[all]{xy}
\usepackage{extarrows}
\usepackage{tabu}
\usepackage{tikz-cd}

\usepackage{txfonts, pxfonts}

\usepackage{amsthm}
\usepackage{amsmath}

\newtheorem{thm}{Theorem}[section]

\newtheorem{alphthm}{Theorem}			

\newtheorem{alphcor}[alphthm]{Corollary}

 \theoremstyle{definition}
  \newtheorem{defn}[thm]{Definition}

 \theoremstyle{remark}

\newtheorem*{claim*}{Claim}

\def\NN{\mathbb N}

\DeclareMathOperator{\supp}{supp}
\newcommand{\norm}[1]{\|#1\|}

\begin{document}

\title{A weight-free characterisation of Yu's Property A}

\author{Jiawen Zhang and Jingming Zhu}

\address[Jiawen Zhang]{School of Mathematical Sciences, Fudan University, 220 Handan Road, Shanghai, 200433, China.}
\email{jiawenzhang@fudan.edu.cn}

\address[Jingming Zhu]{School of Data Sciences, Jiaxing University, 899 Guangqiong Road, Jiaxing, 314000, China.}
\email{jingmingzhu@zjxu.edu.cn}

\date{}
\subjclass[2010]{}

\thanks{JZ was partly supported by NSFC 12422107, and the National Key R{\&}D Program of China 2022YFA100700.}

\baselineskip=16pt

\begin{abstract}
In this short note, we give a complete answer to the question of when the generalised F\o lner sets exhibiting property A can be chosen to be subsets of the space itself. More precisely, we prove that this holds for any discrete metric space of bounded geometry.
\end{abstract}

\maketitle

\section{Introduction}

The notion of Property A is a coarse geometric property introduced by Yu in his groundbreaking work \cite{Yu00}, as a coarse analogue to the notion of amenability. Yu showed that Property A implies the coarse Baum-Connes conjecture, which is a central conjecture in higher index theory and provides an algorithm to compute the $K$-theories of Roe algebras. Roe algebras were introduced by Roe, and their $K$-theories serve as receptacles for higher indices of elliptic differential operators on open manifolds (see \cite{Roe88}). This innovative perspective led to fruitful applications in geometry and topology (see, \emph{e.g.}, \cite{Roe96, WY20}).


Recall that a discrete metric space $(X,d)$ has \emph{property A} if for any $R,\varepsilon>0$ there exists a family $\{A_x\}_{x\in X}$ of finite, non-empty subsets of $X \times \NN$ satisfying:
\begin{itemize}
  \item for any $x,y\in X$ with $d(x,y) \leq R$, we have
  $\frac{|A_x \triangle A_y|}{|A_x \cap A_y|} < \varepsilon$;
  \item there exists an $S>0$ such that for any $x\in X$ and $(y,n)\in A_x$, then $d(x,y)\leq S$.
\end{itemize}
Note that the sets $A_x$ above are allowed to be taken in $X \times \NN$, rather than just the space $X$ itself. This plays a key role in allowing for non-injective maps between coarse spaces. In particular in the proof that property A is inherited by subspaces (see \cite{Wil09b, Yu00}), one needs to be able to retract neighbourhoods onto the subspace.


Let us compare with the definition of amenability: a countable discrete group $G$ is \emph{amenable} if for any $R,\varepsilon>0$ there exists a finite, non-empty subset $A$ of $G$ such that for any $x,y\in G$ with $d(x,y) \leq R$, we have $\frac{|xA \triangle yA|}{|xA \cap yA|} < \varepsilon$. It is clear that the translations $A_x:=\{xg: g\in A\}$ provide the required property A sets. While in this case, the sets $A_x$ can be taken in the group itself. 

This suggests the following simpler version of property A:

\begin{defn}[{\cite[Definition 1.1]{NWZ24}}]\label{defn:naive property A}
A discrete metric space $(X,d)$ has \emph{naive property A} if for any $R,\varepsilon>0$ there exists a family $\{A_x\}_{x\in X}$ of finite and non-empty subsets of $X$ satisfying:
\begin{enumerate}
  \item[(a)] for any $x,y\in X$ with $d(x,y) \leq R$, we have
  $\frac{|A_x \triangle A_y|}{|A_x \cap A_y|} < \varepsilon$;
  \item[(b)] there exists an $S>0$ such that for any $x\in X$ and $y\in A_x$, then $d(x,y)\leq S$.
\end{enumerate}
\end{defn}

Note that the only difference between Definition \ref{defn:naive property A} and property A is that the family of sets $\{A_x\}_{x\in X}$ are required to be in $X$ rather than $X \times \NN$. It is obvious that naive property A implies property A, however, the converse is unclear.

%


In a recent work \cite{NWZ24}, Niblo, Wright and the first-named author proved the converse for a large class of spaces. To state their result, recall that for a metric space $(X,d)$, the \emph{$r$-Rips complex of $X$} is the simplicial complex whose simplices are $\{x_0,x_1, \cdots,x_n\} \subseteq X$ with $d(x_i, x_j) \leq r$ for any $i,j=0,1,\cdots,n$. A metric space is said to \emph{coarsely have unbounded components} if for some $r>0$, every connected component of the $r$-Rips complex is unbounded. Also $(X,d)$ is said to have \emph{bounded geometry} if $\sup_{x\in X} |B_X(x,R)|$ is finite for any $R>0$, where $B_X(x,R):=\{y\in X: d(y,x) \leq R\}$.

It was proved in \cite{NWZ24} that property A is the same as naive property A for the following classes of discrete metric spaces of bounded geometry:
\begin{itemize}
 \item metric spaces which coarsely have unbounded components;
 \item countable groups with proper word length metric;
 \item box spaces of residually finite groups.
\end{itemize}
However, there are still some classes of metric spaces which cannot be covered by the results in \cite{NWZ24}, \emph{e.g.}, coarse disjoint unions of general finite metric spaces. 


%


In this paper, we provide the following complete answer, which makes up the missing piece to the relation of property A and  naive property A:

\begin{alphthm}\label{thm:naive A=A}
Let $X$ be a discrete bounded geometry metric space. Then $X$ has property A \emph{if and only if} $X$ has naive property A.
\end{alphthm}

Our proof partially makes use of the technique in \cite{NWZ24} to construct flows in the space to redistribute mass out towards infinity, which fits well for unbounded components in the Rips complex. While for bounded components (which might \emph{not} have a uniform bound of diameter), we develop a new approach by adding an infinite tail to each bounded component and tailoring the naive property A sets back to the original space. 

Consequently, we obtain the following permanence properties:

\begin{alphcor}
Naive property A is preserved under taking subspaces, finite unions, finite Cartesian products and coarse equivalences.
\end{alphcor}

Let us fix some notions for later use. For a metric space $(X,d)$ and $R>0$, a subspace $A \subseteq X$ is called \emph{$R$-connected} if it is contained in a connected component of the $R$-Rips complex of $X$. An \emph{$R$-connected component} of $X$ is a connected component in the $R$-Rips complex of $X$. By an \emph{$\NN$-valued $0$-chain} on $X$, we mean a map $a: X \to \NN$. Denote its $\ell^1$-norm by $\|a\|_1:=\sum_{x\in X} a(x)$. 

\medskip

\noindent \textbf{Acknowledgement.} The first named author would like to thank Graham Niblo and Nick Wright for helpful comments and continuous support.

\section{Proof of Theorem \ref{thm:naive A=A}}

The whole section is devoted to the proof of Theorem \ref{thm:naive A=A}. It suffices to show that property A implies naive property A, and we divide the proof into several steps.

\medskip

\noindent\textbf{Step I. Reduction from Property A.}

Given $R>0$ and $\varepsilon>0$, $X$ having property A implies that there exists a family of finite, non-empty subsets $\{A_x\}_{x\in X}$ of $X\times \NN$ as specified in the definition of property A. By Theorem 1.2.4 $``(3)\Rightarrow(1)''$ in \cite{Wil09b}, we can find an $M\in \NN$ such that the family $\{A_x\}_{x\in X}$ can be chosen to be subsets of $X \times \{0,1,2,\ldots, M-1\}$. Let us fix $R,\varepsilon$ and $M$ from now on, and our aim is to construct a family of finite, non-empty subsets $\{\widetilde{A}_x\}_{x\in X}$ of $X$ satisfying the conditions in Definition \ref{defn:naive property A} for $R,\varepsilon$.

For each $x\in X$, define an $\NN$-valued $0$-chain $a_x$ by $a_x(z)=|A_x\cap (\{z\}\times \NN)|$ for $z\in X$. By the argument in the proof of \cite[Theorem A]{NWZ24}, the conditions on the sets $A_x$ translate as:
\begin{enumerate}[(i)]
  \item for any $x,y\in X$ with $d(x,y) \leq R$ we have
  \[\dfrac{\norm{a_x-a_y}_1}{\norm{a_x \land a_y}_1} < \varepsilon\]
  where $a\land b$ denotes the pointwise minimum of two chains;
  \item there exists an integer $S>0$ such that for any $x\in X$, $\supp (a_x) \subseteq B_X(x,S)$;
  \item there exists an integer $L>0$ (depending on $M$ and $S$) such that for any $x\in X$, $\norm{a_x}_1 < L$. 
\end{enumerate}
Without loss of generality, we assume $S>R$. Fix such $S$ and $L$.

\medskip

\noindent\textbf{Step II. Constructing a new metric space $(\widetilde{X},\tilde{d})$.}

Let us decompose $X$ into its $S$-connected components, \emph{i.e.}, write:
\[
X=\bigsqcup_{\lambda\in\Lambda} X_{\lambda},
\]
where each $X_\lambda$ is $S$-connected and $d(X_\lambda, X_{\lambda'}) > S$ for $\lambda\neq \lambda'$ in $\Lambda$. According to (ii) in Step I, we have $\supp a_x \subseteq X_\lambda$ for $x\in X_\lambda$. Denote:
\[
  \Lambda_0:=\left\{\lambda\in\Lambda: X_{\lambda} \text{ is unbounded}\right\} \quad \text{and} \quad \Lambda_1:=\left\{\lambda\in\Lambda: X_{\lambda} \text{ is bounded}\right\}.
\]
To use the techniques in \cite{NWZ24}, we need ``add an infinite tail'' to each bounded component. To achieve, firstly choose a point $x_\lambda \in X_\lambda$ for each $\lambda \in \Lambda_1$. Denote $Y_\lambda:=\{\lambda\} \times \NN$, and write $y_\lambda^{(n)}:=(\lambda,n)$ for $n\in \NN$. Denote
\[
\widetilde{X}_{\lambda}:= (X_{\lambda} \sqcup Y_{\lambda})/\sim, \quad \text{where} \quad x_\lambda \sim y_\lambda^{(0)}.
\]
Intuitively, $\widetilde{X}_{\lambda}$ comes from $X_\lambda$ by adding an infinite sequence attached to $x_\lambda$. Now define a metric $d_\lambda$ on $\widetilde{X}_{\lambda}$ as follows: for any $x,x' \in \widetilde{X}_{\lambda}$, set
\begin{equation*}
    d_{\lambda}(x,x') := 
    \begin{cases}
    ~d(x,x'),  & \text{if} \quad x,x'\in X_{\lambda};\\
    ~d(x,x_{\lambda})+j\cdot S, & \text{if} \quad x\in X_{\lambda}~\text{and}~x'=y_{\lambda}^{(j)}~\text{for some}~j\in\NN;\\
    ~d(x_{\lambda}, x')+j\cdot S, & \text{if} \quad x'\in X_{\lambda}~\text{and}~x=y_{\lambda}^{(j)}~\text{for some}~j\in\NN;\\
    ~|j-i|\cdot S, & \text{if} \quad x=y_{\lambda}^{(i)} \text{ and }x'=y_{\lambda}^{(j)} \text{ for some } i,j\in\NN.
    \end{cases}
\end{equation*}    
It is easy to check that $d_{\lambda}$ is a well-defined metric on $\widetilde{X}_{\lambda}$. 
For convenience, also set $\widetilde{X}_{\lambda}=X_{\lambda}$ and $d_{\lambda}:=d|_{\widetilde{X}_{\lambda}}$ for any $\lambda\in\Lambda_0$.

Define a new metric space $\widetilde{X}$ as follows: Denote
\begin{equation}\label{EQ:construction of new space}
\widetilde{X}:= \bigsqcup_{\lambda\in\Lambda}  \widetilde{X}_{\lambda}
\end{equation}
equipped with a metric $\tilde{d}$ such that for any $x, x' \in \widetilde{X}$, we set
\begin{small}
\begin{equation*}
    \tilde{d}(x,x') := 
    \begin{cases}
    ~d(x,x'),  & \text{if} \quad x,x'\in X;\\
    ~d(x,x_{\lambda'})+d_{{\lambda'}}(x_{\lambda'},x'), & \text{if} \quad x\in X_{\lambda}, x'\in Y_{\lambda'} \text{ for some } \lambda\in \Lambda, \lambda'\in \Lambda_1;\\
     ~d_{\lambda}(x,x_{\lambda}) + d(x_{\lambda},x'), & \text{if} \quad x\in Y_{\lambda}, x'\in X_{\lambda'} \text{ for some } \lambda\in \Lambda_1, \lambda'\in \Lambda;\\  
    ~d_{{\lambda}}(x,x_{\lambda})+d(x_{\lambda},x_{\lambda'})+d_{{\lambda'}}(x_{\lambda'},x'), & \text{if} \quad x\in Y_{\lambda}, x'\in Y_{\lambda'} \text{ for some }\lambda \neq \lambda'\in \Lambda_1;\\
    ~d_\lambda(x,x'), & \text{if} \quad x, x' \in Y_{\lambda} \text{ for some }\lambda\in \Lambda_1.
    \end{cases}
\end{equation*}    
\end{small}
It is clear that $\tilde{d}$ is a well-defined metric on $\widetilde{X}$ which also has bounded geometry, and $\tilde{d}|_{X_\lambda} = d_\lambda$ on $\widetilde{X}_{\lambda}$ for each $\lambda \in \Lambda$. Moreover, each $\widetilde{X}_{\lambda}$ is an $S$-connected component of $\widetilde{X}$ and hence, Equation (\ref{EQ:construction of new space}) decomposes $\widetilde{X}$ into its $S$-connected components. By construction, each $S$-connected component $\widetilde{X}_{\lambda}$ is unbounded, \emph{i.e.}, $(\widetilde{X},\tilde{d})$ coarsely has unbounded components defined in \cite[Definition 1.2]{NWZ24}.

\medskip

\noindent\textbf{Step III. Flow on $(\widetilde{X},\tilde{d})$.}

Now we use the techniques in \cite[Section 4]{NWZ24} for the space $(\widetilde{X}, \tilde{d})$. For each $\widetilde{X}_\lambda$, we choose a maximal tree in its $S$-Rips complex, an embedded ray and a point at infinity. Orient the edges of the tree to point towards the chosen basepoint. Then for every $x \in \widetilde{X}$ there is exactly one edge leaving $x$, whose other end we denote $\sigma(x)$. This yields a map $\sigma: \widetilde{X} \to \widetilde{X}$ such that for every $x \in \widetilde{X}$ the sequence $x,\sigma(x),\sigma(\sigma(x)),\cdots$ is proper, \emph{i.e.}, any bounded set contains only finitely many of these points.

Given an $\NN$-valued $0$-chain $a$ on $\widetilde{X}$, let $b(a)$ denote the characteristic function of the support of $a$, and $t(a):=a-b(a)$. Define an $\NN$-valued $0$-chain $s_1(a)$ by:
\[
s_1(a)(x)=b(a)(x)+\sum_{y\in \widetilde{X}: \sigma(y)=x} t(a)(y) \quad \text{for} \quad x\in \widetilde{X}.
\]
Define $s_n$ to be the $n$-fold composition of $s_1$. 
By \cite[Claim 4.2]{NWZ24}, $s_n(a)$ is $\{0,1\}$-valued for $n\geq \norm{a}_1 \cdot \norm{t(a)}_1$. 
Write $s_\infty(a)$ for $s_n(a)$ when $n\geq  \norm{a}_1 \cdot \norm{t(a)}_1$. 
Since the map $s_1$ preserves the $\ell^1$ norm, then $|\supp (s_\infty(a))|=\|a\|_1$.

Now we apply $s_\infty$ to the family $\{a_x\}_{x\in X}$ in Step I. For $x\in X$, take $\lambda \in \Lambda$ such that $x\in X_\lambda$. We have $|\supp (s_\infty(a_x))| = \|a_x\|_1 < L$ by condition (iii) in Step I, and hence $s_\infty(a_x) = s_{L^2}(a_x)$. 
Note that the support of $s_1(a_x)$ lies in the $S$-neighbourhood of the support of $a_x$. Combining with the fact that $\supp(a_x) \subseteq X_\lambda$ and $\widetilde{X}_\lambda$ is an $S$-connected component, we then have $\supp (s_\infty(a_x)) \subseteq \widetilde{X}_\lambda$. Moreover, according to condition (ii) in Step I, we have
\begin{equation}\label{EQ:bdd}
\supp (s_\infty(a_x))\subseteq B_{\widetilde{X}_\lambda}(x, S+S\cdot L^2).
\end{equation}
Take $N:=L^2+2 \in \NN$.
Then for any $x\in X_\lambda$ with $\lambda\in\Lambda_1$, we have
\begin{equation}\label{EQ:restriction}
\supp (s_\infty(a_x))\cap \left\{y_{\lambda}^{(N+1)},y_{\lambda}^{(N+2)},\cdots\right\}=\emptyset.
\end{equation}

On the other hand, applying the same argument for \cite[Claim 4.3]{NWZ24} and thereafter, we obtain
\begin{equation}\label{esti}
\dfrac{|\supp(s_\infty(a_x)) ~\triangle~ \supp(s_\infty(a_y))|}{|\supp(s_\infty(a_x)) \cap \supp(s_\infty(a_y))|} = \dfrac{\norm{s_\infty(a_x)-s_\infty(a_y)}_1}{\norm{s_\infty(a_x) \land s_\infty(a_y)}_1} \leq \dfrac{\norm{a_x-a_y}_1}{\norm{a_x \land a_y}_1} < \varepsilon
\end{equation}
for any $x,y\in X$ with $d(x,y) \leq R$ thanks to condition (i) in Step I.


\medskip

\noindent\textbf{Step IV. Tailoring $\supp(s_\infty(a_x))$ back to $X$.}

Set $\Lambda_{1,+}:= \{\lambda\in \Lambda_1: X_{\lambda} \setminus B_{X_\lambda}(x_{\lambda}, 3S+4S\cdot N)\neq \emptyset\}$. For $\lambda \in \Lambda_{1,+}$, $X_\lambda$ being $S$-connected implies that there exist mutually different points $z_{\lambda}^{(1)},z_{\lambda}^{(2)},\cdots,z_{\lambda}^{(N)}$ in $X_{\lambda}$ such that
\[
\left\{z_{\lambda}^{(1)},z_{\lambda}^{(2)},\cdots,z_{\lambda}^{(N)}\right\}\subset B_{X_\lambda}(x_{\lambda}, 3S+4S\cdot N)\setminus B_{X_\lambda}(x_{\lambda}, 3S+3S\cdot N).
\]
For $\lambda \in \Lambda_1$, define a subspace $\widetilde{Z}_{\lambda} \subseteq \widetilde{X}_{\lambda}$ by
\[
\widetilde{Z}_{\lambda}:= \left(X_{\lambda} \sqcup \left\{y_\lambda^{(0)}, y_\lambda^{(1)}, \cdots, y_\lambda^{(N)}\right\}\right)\big/\sim.
\]
For convenience, also set $\widetilde{Z}_{\lambda}:=X_\lambda$ for $\lambda \in \Lambda_0$. Then Equality (\ref{EQ:restriction}) shows that if $x\in X_\lambda$, then $\supp(s_\infty(a_x)) \subseteq \widetilde{Z}_{\lambda}$.

Given $\lambda \in \Lambda$ and non-empty $A_\lambda \subseteq \widetilde{Z}_{\lambda}$, we define
\begin{equation*}
    \Phi_\lambda(A_\lambda) := 
    \begin{cases}
    ~A_\lambda,  & \text{if} \quad \lambda \in \Lambda_0;\\
    ~X_\lambda, & \text{if} \quad \lambda \in \Lambda_1\setminus \Lambda_{1,+};\\
    ~(A_\lambda \cap X_\lambda) \cup \left\{z_{\lambda}^{(i)}: i \in \NN\setminus \{0\} \text{ with } y_{\lambda}^{(i)} \in A_\lambda\right\}, & \text{if} \quad \lambda \in \Lambda_{1,+}.
    \end{cases}
\end{equation*}  
Then for any $x\in X_\lambda$, we define
\[
\widetilde{A}_x:= \Phi_\lambda(\supp(s_{\infty}(a_x))).
\]
It remains to verify that the family $\{\widetilde{A}_x\}_{x\in X}$ satisfies condition (a) and (b) in Definition \ref{defn:naive property A} for $R$ and $\varepsilon$.

\medskip

\noindent\textbf{Step V. Verifying $\{\widetilde{A}_x\}_{x\in X}$ provides naive property A sets.}

Fix $x,y \in X$ with $d(x,y) \leq R$. Take $\lambda \in \Lambda$ such that $x \in X_\lambda$, which implies that $y\in X_\lambda$ as well. We divide into three cases.

\noindent \emph{Case 1}. $\lambda \in \Lambda_0$. Then $\widetilde{X}_\lambda = X_\lambda$ and $\widetilde{A}_x= \supp(s_{\infty}(a_x))$. Hence condition (a) and (b) in Definition \ref{defn:naive property A} holds thanks to Equality (\ref{EQ:bdd}) and (\ref{esti}).

\noindent \emph{Case 2}. $\lambda \in \Lambda_1 \setminus \Lambda_{1,+}$. Then $\widetilde{A}_x= \widetilde{A}_y = X_\lambda =B_{X_\lambda}(x_\lambda, 3S+4S\cdot N) \subseteq B_{X}(x, 6S+8S\cdot N)$. Hence condition (a) and (b) in Definition \ref{defn:naive property A} holds as well.

\noindent \emph{Case 3}. $\lambda \in \Lambda_{1,+}$. Firstly, we check condition (b). If $\supp(s_{\infty}(a_x)) \subseteq X_\lambda$, then $\Phi_\lambda(\supp(s_{\infty}(a_x))) = \supp(s_{\infty}(a_x))$.
Otherwise, there exists $i=1,2,\cdots, N$ such that $y_\lambda^{(i)} \in \supp(s_{\infty}(a_x))$. Then we have
\begin{align*}
\widetilde{A}_x &= \Phi_\lambda(\supp(s_{\infty}(a_x))) \subseteq \left( \supp(s_{\infty}(a_x)) \cap X_\lambda \right) \cup \left\{z_\lambda^{(0)}, \cdots, z_\lambda^{(N)}\right\} \\
& \subseteq \left( \supp(s_{\infty}(a_x)) \cap X_\lambda \right) \cup \left(B_{\widetilde{X}_\lambda}(y_\lambda^{(i)}, 3S+5NS) \cap X_\lambda\right).
\end{align*}
Hence $\widetilde{A}_x$ is contained in the $(3S+5NS)$-neighbourhood of $\supp(s_{\infty}(a_x))$. According to Equality (\ref{EQ:bdd}), we have
\[
\widetilde{A}_x \subseteq B_{\widetilde{X}_\lambda}(x,S+SL^2 + 3S+5NS) \cap X \subseteq B_X(x,4S+6NS),
\]
which concludes condition (b) in Definition \ref{defn:naive property A}.

It remains to check condition (a). We divide further into two cases.

\noindent \emph{Case 3a}. Both $\supp(s_{\infty}(a_x)) \cap \left\{y_\lambda^{(1)}, \cdots, y_\lambda^{(N)}\right\} = \emptyset$ and $\supp(s_{\infty}(a_y)) \cap \left\{y_\lambda^{(1)}, \cdots, y_\lambda^{(N)}\right\} = \emptyset$.\\
In this case, we have
\[
\widetilde{A}_x = \Phi_\lambda(\supp(s_{\infty}(a_x))) = \supp(s_{\infty}(a_x))
\]
and
\[
\widetilde{A}_y = \Phi_\lambda(\supp(s_{\infty}(a_y))) = \supp(s_{\infty}(a_y)).
\]
Combining with Equality (\ref{esti}), we obtain
\[
\dfrac{\left|\widetilde{A}_x ~\triangle~ \widetilde{A}_{y}\right|}{\left|\widetilde{A}_x \cap \widetilde{A}_{y}\right|} = \dfrac{|\supp(s_\infty(a_x)) ~\triangle~ \supp(s_\infty(a_y))|}{|\supp(s_\infty(a_x)) \cap \supp(s_\infty(a_y))|}  < \varepsilon.
\]

\noindent \emph{Case 3b}. Either $\supp(s_{\infty}(a_x)) \cap \left\{y_\lambda^{(1)}, \cdots, y_\lambda^{(N)}\right\} \neq \emptyset$ or $\supp(s_{\infty}(a_y)) \cap \left\{y_\lambda^{(1)}, \cdots, y_\lambda^{(N)}\right\} \neq \emptyset$. \\
Without loss of generality, assume $\supp(s_{\infty}(a_x)) \cap \left\{y_\lambda^{(1)}, \cdots, y_\lambda^{(N)}\right\} \neq \emptyset$. By Equality (\ref{EQ:bdd}), we have $B_{\widetilde{X}_\lambda}(x,S+SL^2) \cap \left\{y_\lambda^{(1)}, \cdots, y_\lambda^{(N)}\right\} \neq \emptyset$. Assume $y_\lambda^{(i)} \in B_{\widetilde{X}_\lambda}(x,S+SL^2)$ for some $i=1,2,\cdots,N$. Then
\begin{equation}\label{EQ:x and xlambda}
d(x,x_\lambda) \leq d_\lambda(x,y_\lambda^{(i)}) + d_\lambda(y_\lambda^{(i)},x_\lambda) \leq S+SL^2 + NS \leq S+2NS.
\end{equation}
Hence we obtain
\[
B_{\widetilde{X}_\lambda}(x,S+SL^2) \cap X_{\lambda} \subseteq B_{X_{\lambda}}(x_\lambda, S+SL^2 + S+2NS) \subseteq B_{X_\lambda}(x_\lambda, 2S+ 3NS),
\]
which shows that
\[
\left(\supp(s_{\infty}(a_x)) \cap X_\lambda \right) \cap \left\{z_{\lambda}^{(1)},z_{\lambda}^{(2)},\cdots,z_{\lambda}^{(N)}\right\} = \emptyset.
\]
Since $d(x,y) \leq R <S$, we have $\supp (s_\infty(a_y)) \subseteq X_\lambda$ as well. Moreover, using Equality (\ref{EQ:bdd}) and (\ref{EQ:x and xlambda}), we obtain:
\begin{align*}
\supp (s_\infty(a_y)) &\subseteq B_{\widetilde{X}_\lambda}(y, S+SL^2) \subseteq B_{\widetilde{X}_\lambda}(x, 2S+SL^2)
 \subseteq B_{\widetilde{X}_\lambda}(x_\lambda, 2S+SL^2 + S+2NS)\\
& \subseteq B_{\widetilde{X}_\lambda}(x_\lambda, 3S+3NS).
\end{align*}
Hence we have
\[
\left(\supp(s_{\infty}(a_y)) \cap X_{\lambda}\right) \cap \left\{z_{\lambda}^{(1)},z_{\lambda}^{(2)},\cdots,z_{\lambda}^{(N)}\right\} = \emptyset.
\]
Therefore, we obtain
\begin{align*}
\left|\widetilde{A}_x ~\triangle~ \widetilde{A}_y\right| &=  \left|\Phi_\lambda(\supp(s_{\infty}(a_x))) ~\triangle~ \Phi_\lambda(\supp(s_{\infty}(a_y)))\right| \\
&=\left|(\supp(s_{\infty}(a_x)) \cap X_{\lambda}) ~\triangle~ (\supp(s_{\infty}(a_y)) \cap X_{\lambda})\right| \nonumber \\
& \quad \quad + \left|\left\{z_\lambda^{(i)}: i \in \NN \setminus \{0\} \text{ with } y_\lambda^{(i)} \in \supp(s_{\infty}(a_x)) ~\triangle~ (\supp(s_{\infty}(a_y)) \right\}\right| \nonumber\\ 
&= \left|\left(\supp(s_{\infty}(a_x)) ~\triangle~ (\supp(s_{\infty}(a_y))\right) \cap X_\lambda\right| \nonumber\\
& \quad \quad + \left|\left(\supp(s_{\infty}(a_x)) ~\triangle~ (\supp(s_{\infty}(a_y))\right) \cap \left\{y_\lambda^{(1)}, \cdots, y_\lambda^{(N)}\right\}\right|\nonumber \\
&= \left|\supp(s_{\infty}(a_x)) ~\triangle~ (\supp(s_{\infty}(a_y))\right|, \nonumber
\end{align*}
and 
\begin{align*}
\left|\widetilde{A}_x \cap \widetilde{A}_y\right| &=  \left|\Phi_\lambda(\supp(s_{\infty}(a_x))) \cap \Phi_\lambda(\supp(s_{\infty}(a_y)))\right| \\
&=\left|(\supp(s_{\infty}(a_x)) \cap X_{\lambda}) \cap (\supp(s_{\infty}(a_y)) \cap X_{\lambda})\right| \nonumber \\
& \quad \quad + \left|\left\{z_\lambda^{(i)}: i \in \NN \setminus \{0\} \text{ with } y_\lambda^{(i)} \in \supp(s_{\infty}(a_x)) \cap (\supp(s_{\infty}(a_y)) \right\}\right| \nonumber\\ 
&= \left|\left(\supp(s_{\infty}(a_x)) \cap (\supp(s_{\infty}(a_y))\right) \cap X_\lambda\right| \nonumber\\
& \quad \quad + \left|\left(\supp(s_{\infty}(a_x)) \cap (\supp(s_{\infty}(a_y))\right) \cap \left\{y_\lambda^{(1)}, \cdots, y_\lambda^{(N)}\right\}\right|\nonumber \\
&= \left|\supp(s_{\infty}(a_x)) \cap (\supp(s_{\infty}(a_y))\right|. \nonumber
\end{align*}
Combining the above together, we obtain
\[
\dfrac{\left|\widetilde{A}_x ~\triangle~ \widetilde{A}_{y}\right|}{\left|\widetilde{A}_x \cap \widetilde{A}_{y}\right|}  = \dfrac{\left|\supp(s_\infty(a_x)) ~\triangle~ \supp(s_\infty(a_y))\right|}{\left|\supp(s_\infty(a_x)) \cap \supp(s_\infty(a_y))\right|}  < \varepsilon.
\]
This finishes the proof.

\bibliographystyle{plain}
\bibliography{NaiveA}

\end{document}